\documentclass[a4paper, 11pt]{article}
\usepackage{graphicx}

\begin{document}

\title{Equal confidence weighted expectation value estimates}
\author{Fetze Pijlman\\ \begin{small}fetze.pijlman@philips.com\end{small} \\ \small Philips Lighting Research\\ \small High Tech Campus 7, Eindhoven, The Netherlands}

\maketitle

\begin{abstract}
In this article the issues are discussed with the Bayesian approach, least-square fits, and most-likely
fits. Trying to counter these issues, a 
method, based on weighted confidence, is proposed for estimating probabilities and other observables.
This method sums over different model parameter combinations but does not require the need for making assumptions on priors or underlying probability functions. Moreover, by construction
the results are invariant under reparametrization of the model parameters.
In one case the result appears similar
as in Bayesian statistics but in general there is no agreement. 
The binomial distribution is also studied which turns out to be useful for making predictions on 
production processes without the need to make further assumptions. In the last part, the case
of a simple linear fit (a multi-variate example) is studied using the standard approaches and the confidence weighted approach.

\end{abstract}

\section{Introduction}

In statistics we like to analyze data combined with some additional information in order to make predictions on the underlying population and/or future expectations. Consider the following example.
When comparing an assumed probability function $p(x|\tau)$, where $\tau$ is a continuous parameter, to experimental data it is common practice in classical statistics
to define a certain confidence interval (how certain does one wants to be) after which interval ranges 
for
the parameter $\tau$ can be ruled out. One remains with possible values for $\tau$ (see for example Ref.~\cite{prosper}). Although it is certainly nice that we can rule out
certain ranges for $\tau$, this by itself is not sufficient to give an expectation value for some observable, for example, the mean of the population.
For an introduction into statistics and probability theory the reader is referred to Ref.~\cite{larsen}.

In order to come with an estimate for the mean one can use another common
technique which is making a least-square fit. Although least-square fits enables one to easily derive some value for $\tau$, care should be taken with its interpretation.
Least-square fits 
happen in some cases to correspond to the most likely scenario for the 
probability function (see for example Ref.~\cite{least_square}) and although most-likely certainly contains valuable information, it is often used to make further 
predictions which can easily fail as shown by the following example.

Consider for example 10 vases where the first vase has 91 red balls and 9 blue balls, 
the second has 92 red balles and 8 red balles, etc.,
and the last vase has 100 red balls. Blindfolded a ball is taken from a vase which turns out to be red. The most likely
vase from which this red ball was taken is the vase having the 100 red balls. 
If we use this most-likely scenario/fit for making further 
predictions we could easily be making a huge mistake. Suppose we would estimate the probability for
drawing a blue ball from the same vase from which we drew the red ball we would estimate this 
probability to be zero. However, although the most-likely probability after having drawn the red ball is indeed
zero, the expecation value is certainly nonzero as there is a nonzero chance that we could have drawn
the balls from one of the other vases.

Bayesian statistics is one attractive option
for addressing the above problem. One makes an inventory of 
possible hypotheses and uses the data to calculate
the probability that a certain hypothesis is true. 
Using all those probabilities one can derive a kind of 
expected probability function. The basic Equations are the following (see e.g. Ref~\cite{prosper}):
\begin{equation}
p(x_2 | x) = \int \rm{d}\tau\ p(x_2| \tau)\  \mathcal{P}(\tau | x), 
\label{invariance}
\end{equation}
where
\begin{equation}
P(\tau | x) \propto \int {\rm d}x\ p(x| \tau)\ \mathcal{P}(\tau).
\end{equation}
$\mathcal{P}(\tau)$ which describes the initial probability for a certain 
hypothesis before any data was taken and is known as the prior. 

One desired property is the so-called invariance under reparametrization. Consider for example
the two Equations above in the case of a normal distribution with a known mean and unknown standard
deviation $\sigma$. One experimentator could make the choice for considering $p(x| \sigma)$ and integrate in 
Equation~\ref{invariance}
over $\sigma$ while another experimentator could make another choice such as $p(x| s=\sigma^2)$ and
integrate over $s$. It is desired that both experimenters obtain the same result which can be 
imposed through choosing the prior with certain properties.
There are several popular choices such as 
Jeffreys prior~\cite{jeffrey1, jeffrey2}. Although they may fulfill some necessary requirements such as invariance under reparametrization,
they do lack an ab initio derivation.

There seems to be no fundamental method for estimating an expectation probability function from a set of 
experimental data. When the experimentator assumes a probability function $p(x| \tau)$ both classical as well as Bayesian statistics cannot ab initio provide an expected probability function in case
no other information is known. It is the goal of this article to obtain some estimates for expectation values without making additional assumptions.

\section{Confidence level expectation estimated probabilities}

Although there are statistical approaches with different view points there seems to be a common view 
when
something is probable or not. Suppose that we have a data point at x=0 and we know that it is drawn from a
normal distribution with $\sigma=1$ but with unknown $\mu$. At different confidence levels we can deduct an
interval for the parameter $\mu$ as shown in Table~\ref{Table_confidence}. The Table illustrates that at full confidence one
cannot limit $\mu$ to a finite interval but also that the most likely value for $\mu=0$.

\begin{table}
\begin{center}
\begin{tabular}{rrr}
confidence level & confidence level interval & interval difference wrt previous\\
\hline
10\% & [-0.13, 0.13] & [-0.13,0] \& [0, 0.13]\\
20\% & [-0.25, 0.25] & [-0.25, -0.13] \& [0.13, 0.25]\\
30\% & [-0.39, 0.39] &[-0.39, -0.25] \& [0.25, 0.39]\\
40\% & [-0.52, 0.52] &[-0.52, -0.39] \& [0.39, 0.52]\\
50\% & [-0.67, 0.67] &[-0.67, -0.52] \& [0.52, 0.67]\\
60\% & [-0.84, 0.84] &[-0.84, -0.67] \& [0.67, 0.84]\\
70\% & [-1.04, 1.04] &[-1.04, -0.84] \& [0.84, 1.04]\\
80\% & [-1.28, 1.28] &[-1.28, -1.04] \& [1.04, 1.28]\\
90\% & [-1.64, 1.64] &[-1.64, -1.28] \& [1.28, 1.64]\\
100\%& [-$\infty$, $\infty$] &[-$\infty$, -1.64] \& [1.64, $\infty$]
\end{tabular}
\end{center}
\caption{Confidence interval for a normal distribution with $\sigma=1$ and with unknown mean given
a single data point at $x=0$.\label{Table_confidence}}
\end{table}

When we look at Table~\ref{Table_confidence} from a different point of view then we notice 
that each row makes an equal 
contribution to our confidence (first column shows that each row contributes 10\% with respect to previous row). 
As confidence levels/intervals is at the heart of any statistical method it
makes sense to use this for defining an expectation probability function. Expectation values
could be obtained by considering the intervals in the last column of each row in Table~\ref{Table_confidence} as equally important. It is proposed that intervals contributing equally to confidence shall be
taken as equally contributing to the expectation value.

Before defining the expectation value, the following confidence level limit is defined as
\begin{equation}
\alpha(\vec{x},\tau) \equiv \int_{p(\vec{y}|\tau)>p(\vec{x}|\tau)} \prod_{i=1}^n \rm{d}y_i\ p(\vec{y}| \tau),
\label{definition}
\end{equation}
where $p(\vec{x}|\tau)=\prod_i^n\ p(x_i| \tau)$, $x_i$ are the known measurements, and $\tau$
is a model parameter. Please
note that this definition also nicely deals with probability distributions that have multiple maxima
or that have a maximum at its border.

For each value $\alpha$ one can find either no, one, or more values for $\tau$. 
The number of solutions we indicate with 
$N(\vec{x}, \alpha)$. The expected probability distribution based on equal contributing confidence is now defined as
\begin{equation}
\left< \mathcal{O} \right>_c \equiv \frac{1}{K} \int_{N(\vec{x}, \alpha) \neq 0} \rm{d}\alpha \ \frac{1}{N(\vec{x}, \alpha)} 
\sum_{i=1}^{N(\alpha, \vec{x})} \mathcal{O}(\tau_i(\vec{x}, \alpha)),
\end{equation}
where
\begin{equation}
K = \int_{N(\vec{x}, \alpha) \neq 0} \rm{d}\alpha.
\end{equation}
The amount to the expection to which regions of $\tau$ contribute is determined by $\alpha$. For this
reason the expectation probability function is by definition invariant under reparametrization.

\section{Examples involving normal distribution}

\subsection{Normal distribution with unknown $\mu$}

In order to derive the expectation probability function it turns out to be most
convenient for first compute $\alpha$ after which one can use its derivative with respect to $\mu$ to 
rewrite integral over $\alpha$ into an integral over $\mu$. After translating the integration variables
$y_i$ one can write $\alpha$ as
\begin{equation}
\alpha(\vec{x},\mu) = \int_{p(\vec{y}|0) > p(\vec{x}|\mu)} \prod_{i=1}^n \rm{d}y_i\ p(\vec{y}| 0).
\end{equation}
Using the
Jacobian for spherical coordiantes in $n$ dimensions, $2\pi^{n/2}/\Gamma(n/2)\ r^{n-1}$, see also Ref.~\cite{jacobian1}
, one finds
\begin{equation}
\alpha = \frac{1}{\Gamma (n/2)} \left[ \Gamma(n/2, 0) - \Gamma (n/2, - \ln \prod_k e^{-(x_k-\mu)^2/(2\sigma^2)})
\right],
\end{equation}
where $\Gamma$ is the (upper) incomplete Gamma function. When taking the derivative with respect to $\mu$ one 
obtains
\begin{equation}
\partial_\mu \alpha = \frac{-(\sqrt{2\pi}\sigma)^n}{\sqrt{2}\sigma\ \Gamma(n/2)}\ 
\left[ \sum_i^n \frac{(x_i - \mu)^2}{2\sigma^2} \right]^{n/2 - 1/2}\ p(\vec{x}| \mu).
\end{equation}

For the expected probability function, one obtains
\begin{equation}
\left< p(x_2 | \vec{x} ) \right>_c = \frac{1}{K} \int \rm{d}\mu\ p(x_2| \mu)\ P(\mu, \vec{x})\ p(\vec{x}| \mu),
\end{equation}
where
\begin{equation}
P(\mu, \vec{x}) = \frac{(\sqrt{2\pi}\sigma)^n}{2\sqrt{2}\sigma\ \Gamma(n/2)}\ 
\left[ \sum_i^n \frac{(x_i - \mu)^2}{2\sigma^2} \right]^{n/2 - 1/2}.
\end{equation}
In general the factor $K$ is not equal to $1$ as for multiple different data points it is not feasible
to find $\alpha=0$. However, the factor $K$ can easily be computed by using the property that
$\int {\rm d} x_2 \ \left< p(x_2 | \vec{x} ) \right>_c = 1$.

In the case of a single data point ($n=1$) one finds
\begin{equation}
\left< p(x_2 | x_1 ) \right>_c = \int \rm{d}\mu\ p(x_2| \mu)\ P(\mu, x_1)\ p(x_1| \mu),
\end{equation}
where
\begin{equation}
P( \mu | x_1) = 1.
\end{equation}
The result above is similar to what one could obtain in Bayesian statistics but please do note that
in case of multiple data points that $P$ does contain a dependence on the initial data $x_1$.

\subsection{Normal distribution with unknown $\sigma$}

The treatment of this case is quite similar to the previous case. One finds
\begin{equation}
p(x_2 | \vec{x} ) = \frac{1}{K} \int \rm{d}\sigma\ p(x_2| \sigma)\ P(\sigma , \vec{x})\ p(\vec{x}| \sigma),
\end{equation}
where
\begin{equation}
P(\sigma,\vec{x}) = \frac{(\sqrt{2\pi}\sigma)^n}{\Gamma(n/2)}
\left[ \sum_i^n \frac{(x_i-\mu)^2}{2 \sigma^2}\right]^{n/2-1}\ 
\left| \partial_\sigma p(\vec{x}| \sigma) \right|.
\end{equation}

In case of a single data point, $n=1$, one finds
\begin{equation}
p(x_2 | x_1 ) = \int \rm{d}\sigma\ p(x_2| \sigma)\ P(\sigma, x_1)\ p(x_1| \sigma),
\end{equation}
where
\begin{equation}
P(\sigma, x_1) = 2 \frac{|x_1-\mu|}{\sigma}.
\end{equation}
When comparing these results with the Bayesian approach one concludes that for $n=1$
it is not possible to identify a prior that does not depend on the
initial measurements.

\section{Binomial distribution}

In a production process one often encounters products that are produced with a certain spread. This
spread can be binned in certain intervals. Some bins are considered to be accepTable while other
bins are seen as failures. A question that often occurs is whether the production process is good enough. This is usually investigated by considering a subset of the production (batch) on which
a statistical analysis is carried out with the purpose of making a claim on the total production (population).

An often used approach is to argue that (1) the products deviations follow a certain distribution,
e.g. normal distribution, (2) a fit is being made with the batch, and (3) using the fit an estimate
is made for the failures in the population.
In this section we try to estimate the expection for a certain bin while not making any
assumptions on underlying the probability distribution. As the integrals were too difficult to
compute analytically the results below were obtained via computing the integrals numerically.

The binomial distribution is defined as
\begin{equation}
p(k|n,p) = {n\choose k}\ p^k\ (1-p)^{n-k}.
\end{equation}
Equation~\ref{definition} in this case becomes
\begin{equation}
\alpha = \int_{p(l|n,p)> p(k|n,p)} {\rm d} l\ p(l|n,p),
\end{equation}
and for the expected probability function we get
\begin{equation}
\left< p(k_2 | k_1) \right> = \frac{1}{K} \int_{N \neq0} {\rm d}p'\ \frac{{\rm d}\alpha}{{\rm d} p'}
p(k_2|n,p').
\end{equation}

\begin{table}
\begin{tabular}{r|rrrr}
\# counts &  \# counts 
& st. dev.&  \# counts &
st. dev. \\
found in bin&  most likely 
& most likely & expected&
expected\\
\hline
0 &0&0&0.89&0.49\\
1 &1&.42&1.68&0.59\\
2 &2 & 0.44 & 2.46 & 0.62\\
3 &3 & 0.47 & 3.24 & 0.64
\end{tabular}
\caption{Most likely and expected bin values in case in total 10 samples were measured with a  number of counts observed in a 
certain bin (column 1).\label{Table10}}
\end{table}

\begin{table}
\begin{tabular}{r|rrrr}
\# counts &  \# counts 
& st. dev.&  \# counts &
st. dev. \\
found in bin&  most likely 
& most likely & expected&
expected\\
\hline
0 &0&0 	&1.12	&0.57\\
5 &5&	.57	&6.5&0.83\\
10 &10 	& 0.66 	& 12.14 & 0.90\\
15 &15 	& 0.71 	& 17.58 & 0.94
\end{tabular}
\caption{Most likely and expected bin values in case in total 100 samples were measured with a  number of counts observed in a 
certain bin (column 1).\label{Table100}}
\end{table}

The integral above can be calculated numerically and results for $n=10$ and $n=100$ are given in 
Table~\ref{Table10} and Table~\ref{Table100}. From a batch of $n$ samples the most left column indicates the number
that was observed (e.g. number of products with a failure), and the other columns indicate what one
can expect for the population based on most-likely or weighted by confidence.

From the Tables one can see that in case a low number is observed from a large batch ($k<<n$) that the
confidence weighted expectation for failures in the population is slightly higher than what would
be expected from a most-likely analysis.

\section{A multi-variate example: linear fit}

\begin{figure}
\begin{center}
\includegraphics[width=10cm,angle=0]{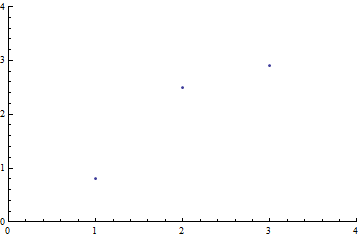}
\caption{An example of data\label{data}. The data points plotted are (1, .8), (2, 2.5), (3, 2.9).}
\end{center}
\end{figure}

Suppose we are confronted with data as given in Figure~\ref{data}. We are told that the data is taken from measurements of which it is 
known that each data point is normally spread ($\sigma=1$) around a linear relation: $y=ax+b$. The question that is posed is
to make a prediction on the quantity $V=(1+|b|)^2$.

When we make a fit based on \emph{least-squares} or one \emph{most-likelyhood} one would find in this case that
$a=1.05$ and $b=-.0333$. The most-likely value for $V$ would be $1.068$. When calculating the standard
error in $V$ via $\sigma_V = |\partial V / /\partial b | \sigma_b$ one would find $\sigma_V = 1.6$.

The interesting aspect here is that the most-likely value for $b$ is close to zero and that any deviation
from zero will result in a larger value for $V$. Any uncertainty in our fit shifts the value for $V$ 
in one particular direction. For this reason one can question whether the most-likely value including
its standard error has some practical relevance.

When we follow the \emph{Bayesian approach} we could assume that the prior $\mathcal{P}(a, b)$
for having certain values for $a$
and $b$ is constant over an interval $[-10,10]$ and zero else where. One would find that
\begin{equation}
\mathcal{P}(a,b | \vec{x} ) \propto \prod_{i=1}^3 P(x_i, y_i) \mathcal{P}(a, b),
\end{equation}
where
\begin{equation}
P(x_i, y_i) = \frac{1}{\sqrt{2\pi}\sigma} e^{-(ax_i+b-y_i)^2/(2)}.
\end{equation}
The expectation value for $V$ can then be computed by
\begin{equation}
\left< V \right> = \int_{-10}^{10} {\rm d}a \int_{-10}^{10}{\rm d}b\ \mathcal{P}(a,b | \vec{x} )\ V(b).
\end{equation}

When carrying out the calculations one finds in this case that the expection value for $V$ becomes
$\left< V \right>=5.77$ with $\sigma = 5.1$. This value is substantially different from the most-likely value which is what we kind of
expected as all deviations would result in a shift to higher values. One problem of this approach
is the assumption of having a constant prior over a finite interval. What is we would have parametrized
our linear relation as $y=a^3 x + b$ or what if we would have chosen the interval $[-20,20]$ instead of the interval $[-10,10]$.
How would our answer change and what would be the correct one?

\begin{figure}
\begin{center}
\includegraphics[width=10cm,angle=0]{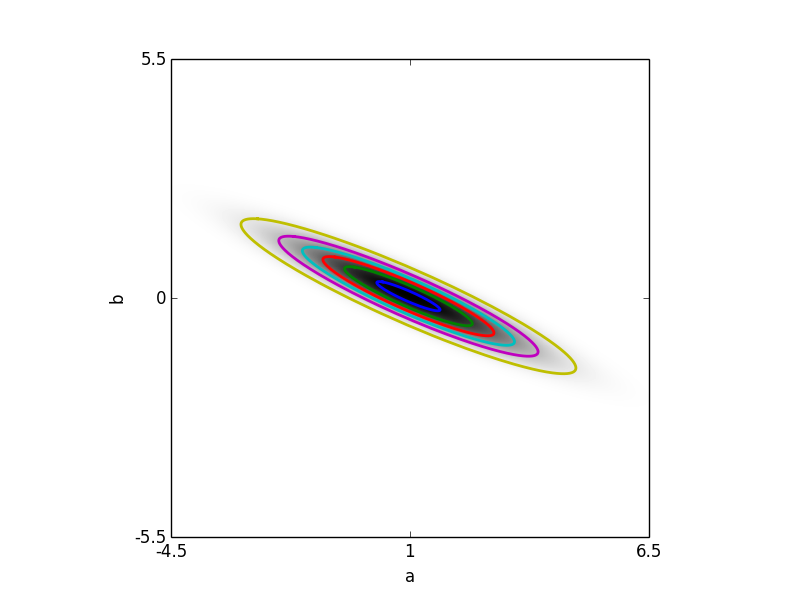}
\caption{A contour plot of $\alpha$ \label{contour}.} 
\end{center}
\end{figure}

Following the suggested approach in this article, weighting with confidence, we can keep the definition for $\alpha$
for this problem (Equation~\ref{definition}). In this case we write $\alpha$ as ($n=3$)
\begin{equation}
\alpha = \frac{1}{\Gamma (n/2)} \left[ \Gamma(n/2, 0) - \Gamma (n/2, - \ln \prod_k e^{-(ax_k+b-y_k)^2/(2)})
\right],
\end{equation}
where $x_k$ and $y_k$ are the measured data points.

In contrast to $\alpha$, the definition for the expectation value
needs to be adapted as the two model parameters (\emph{multi-variate}), $a$ and $b$, can lead to a continuum of solutions
at a certain $\alpha$. This means that $N$ has to be replaced. It seems natural to replace $N$ by 
a line integral in $ab$ space over the contours of equal $\alpha$. We write
\begin{equation}
\left< V \right> = \frac{1}{K} \int {\rm d}\alpha\ \frac{1}{\int_{\alpha} {\rm d}l_{ab}} \ 
\int_{\alpha} {\rm d}l_{ab} \ V(b).
\end{equation}
Please note that the line integral is yet again invariant under reparametrization as
the contour does not change under reparametrization (defined by $\alpha$) and (e.g.)
\begin{eqnarray}
\int_{\alpha} {\rm d}l_{ab} \ f(a,b) &=& \int_{0}^{1} {\rm d}t \sqrt{ (a'(t))^2 + (b'(t))^2 }\ f(a(t),b(t)) \\
&=& \int_0^1 {\rm d}t \sqrt{ ({a}'(h(t)))^2 + ({b}'(h(t)))^2 }\ f(a(h(t)),b(h(t)))\\
&=&\int_{\alpha} {\rm d}l_{a'b'} \ f(a',b').
\end{eqnarray}

The integral above can be computed numerically in the following manner. 
First one determines $\alpha(a,b)$ on a grid.
From this grid one can use scipy.skimage.measure to find the contours at certain levels of
$\alpha$ on the grid (see Figure~\ref{contour}). After having obtained the contours one can carry out the integration, e.g.
using the trapezium rule. Following this method we found $\left<V\right>\approx 5.7$
with $\sigma = 5.1$. This result is very similar to what was achieved with the Bayesian approach.

\section{Conclusions}

In this article the difference between most-likely and expectation has been discussed. Standard
least mean square fits lead in some cases 
to a most likely estimate for the fitted parameters. Although
most-likely is certainly of interest, often we are also interested in what we can expect.

In this article it is proposed to use a kind of summation over hypotheses (or integrating/summing over 
parameters) and taking as weights the amount to which they contribute to the confidence. The whole
idea is that model parameter intervals making an equal contribution to confidence should also
make an equal contribution to the expection value.
It was shown
that the resulting estimate for expectation values does not require a choice for priors. Moreover, the
result is also invariant under reparametrization of the model parameters.

In the case of a normal distribution with a single data point and model parameter the mean, it was found that the result is similar
to Bayesian statistics with a prior equal to $1$. However, in general it was found that the
corresponding prior would depend on the initial data which cannot occur in a Bayesian approach.

The binomial distribution was also studied. This distribution can be of interest
when considering production processes. While not making any assumptions on the underlying
probability distributions, it is nevertheless possible to make some predictions on counts in bins which
differ from most-likely estimates from classical statistics.

In the last section we tried to fit some data with a linear model, $y=ax+b$, having the two parameters.
For this multi-variate case we extended our definition for confidence weighted expectation value. By 
using numerical techniques an expectation value could be obtained which turned out to be close
to a choice that was made for a prior in combination with a Bayesian approach.

\end{document}